# Neue Ergebnisse für das Sammelbilderproblem mit Tauschen und Nachkaufen

Niklas Braband[1], Sonja Braband[1] und Malte Braband[2]

[1]Gymnasium Neue Oberschule, Braunschweig, Deutschland
{niklas.braband|sonja.braband}@no-bs.de
[2]Technische Universität Braunschweig, Deutschland
m.braband@tu-braunschweig.de

**Zusammenfassung.** In diesem Artikel wird das Sammelbilderproblem mit Nachkaufen und Tauschen untersucht. Wir haben mit Hilfe von Simulationen neue Ergebnisse gefunden und sie auch für einen Spezialfall bewiesen. Das Verhältnis der mittleren Anzahl zu kaufender Karten pro Sammler als auch das Verhältnis der Varianz zur Albumgröße hängt in sehr guter Näherung nur vom Anteil der nachkaufbaren Karten ab, aber nicht von der Albumgröße. Damit können Sammler anhand des Prozentsatzes der Nachkaufkarten die mittleren Kosten eines Albums sowie deren Standardabweichung nur mit einer Tabelle und Grundrechenarten bestimmen. Außerdem konnten wir belegen, dass asymptotisch der Effekt des Nachkaufens stärker ist als der des Tauschens.

**Abstract.** This paper focuses on the Coupon Collector's Problem with replacement (limited purchasing of missing stickers) and swapping. We have simulated combined strategies and found new results, which we were able to prove for a particular case. The ratio of the average number of stickers needed to fill the album by the album size (number of different stickers needed) and the number of collectors as well as the ration of the variance and the album size depend approximately only on the percentage of replacement stickers (the limited number of stickers that can be bought directly from the vendor) and the number of collectors, but not on the total album size. Thus collectors can estimate the average cost of completion of an album and its standard deviation just based on basic calculations and a table lookup. Additionally we could show that asymptotocally the effect of replacement is stronger than swapping.

# 1    Einführung

Sammelbilder gibt es zu vielen Themen, z.B. Fussballbilder zur WM. Und es gibt viel Diskussionen über sie, z.B. viele Doppelte, manche scheinen seltener vorzukommen, einige vermuten sogar Betrug der Hersteller. Wir haben begonnen, uns für das Thema zu interessieren [1], als Schweizer Wissenschaftler eine optimale Sammelstrategie [2] vorgeschlagen haben, die auch in der Presse breit veröffentlicht wurde.

Wikipedia [3] schreibt "Das Sammelbilderproblem .... befasst sich mit der Frage, wie viele zufällig ausgewählte Bilder einer Sammelbildserie zu kaufen sind, um eine komplette Bildserie zu erhalten". D.h. es geht darum, eine bestimmte Anzahl B von Bildern vollständig zu sammeln. Man kauft die Bilder normalerweise nicht einzeln, sondern in Päckchen, in denen je P Bilder drin sind. Meistens bietet der Hersteller an, dass man einmalig eine begrenzte Anzahl von Bildern K (zu einem erhöhten Preis) nachkaufen kann, um seine Sammlung zu vervollständigen. Bei Panini war z.B. in Deutschland zur WM 2014 B=640, P=5 und K=50. Weiter ist noch der Preis p eines Päckchens wichtig, sowie der Preis b für die Nachbestellung eines Bildes (im Beispiel p=60 Cent; b=18 Cent).

Weiter ist es noch wichtig, ob man alleine sammelt oder Bilder mit Freunden tauscht. F soll die Anzahl der Freunde sein (einschließlich des Sammlers selbst).

Im klassischen Sammelbilderproblem werden die folgenden Annahmen getroffen:

A1. Die Bilder sind rein zufällig auf die Päckchen verteilt. D.h. der Hersteller mischt bei der Herstellung ordentlich.
A2. In einem Päckchen kommt kein Bild doppelt vor.
A3. Alle Bilder kommen gleich häufig vor. D. h. der Hersteller betrügt nicht durch absichtliche Verknappung von Bildern.
A4. Es wird fair getauscht, d.h. ein Bild gegen ein anderes.
A5. Es gibt keine Rabatte, alle Bilder (außer Nachbestellungen) sind gleich teuer.

Aber man erkennt schnell, dass Sammeln ein teures Vergnügen ist. Denn um das Panini WM-Album zu füllen, hätte man mindestens 128 Päckchen zu 60 Cent kaufen müssen, was 76,80€ gekostet hätte.

Kurz vor der WM haben Sardy und Velenik [2] behauptet, dass die folgende Strategie optimal ist und man im Durchschnitt weniger als 125€ ausgeben muss, um das Album zu füllen:

1.  Kaufe eine Box (bestehend aus 100 Päckchen mit je 5 Bildern)
2.  Kaufe 40 zusätzliche Päckchen und tausche die Doppelten, bis maximal 50

Sticker in der Sammlung fehlen.
3. Kaufe die fehlenden 50 Bilder bei Panini nach.

Eine Box (im Handel auch Display genannt) ist meistens billiger, als die Päckchen einzeln zu kaufen. Dabei geht man in [2] davon aus, dass in einer Box alle Bilder verschieden sind. Sieht man aber im Internet bei Amazon die Kommentare zu den Boxen an, so gibt es Beschwerden über viele doppelte Bilder. Bei den uns von Topps überlassenen Boxen der Serie Match Attax gab es Doppelte. Auch in den Angeboten wird nicht damit geworben, dass alle Bilder in einer Box verschieden sind. Manche meinen, dass es trotzdem weniger Doppelte gäbe als beim Kauf der einzelnen Päckchen. Daher betrachten wir hier die Möglichkeit, Boxen zu kaufen, in dieser Arbeit nicht.

Allerdings konnten die Sardy und Velenik [2] die Kosten für ihre Strategie nur durch Simulation bestimmen. Für jede andere Albumgröße hätte neu simuliert werden müssen. Eine allgemeine Berechnungsformel konnte nicht angegeben werden.

In der Literatur sind viele Ergebnisse über das klassische Sammelbilderproblem [2,3,4,5,6,7] zu finden. Nach [3] gilt für die mittlere Anzahl M zu kaufender Bilder bei einem Sammler $BH(B)$, wobei $H(B)$ die harmonische Summe mit B Summanden bedeutet. D. h. der Faktor f =M/B beträgt hier gerade $H(B)$. Außerdem ist dort angegeben, dass die harmonische Summe gut durch den natürlichen Logarithmus angenähert werden kann, d. h. für große B gilt näherungsweise f=ln(B).

Für das faire Tauschen ist bekannt [4], dass die mittlere Anzahl zu kaufender Bilder bei F Sammlern näherungsweise

$$Bln(B) + B(F-1)ln(ln(B)) \quad (1)$$

beträgt. Dies bedeutet, dass für eine große Anzahl F von Sammlern der Faktor f= M/(BF) gegen ln(ln(B)) strebt.

## 2 Unsere Vermutung

Die Arbeiten aus der Literatur [4,5,6,7] berücksichtigen meistens nur einen Effekt, z. B. Tauschen. Daher haben wir begonnen, mittels Simulationen in R [8] komplexere Versuche durchzuführen, die mehrere Effekte berücksichtigen können.

Dabei haben wir den Effekt des Tauschens und Nachkaufens bei fünf Alben mit sehr unterschiedlicher Größe und zwei Prozentsätzen von Nachkaufkarten K/B ausgewertet. Wir haben zuerst einmal intuitiv 1000 Simulationen für jedes Album und jede

Anzahl Sammler F durchgeführt. Dann haben wir den Faktor f gegen die Anzahl der Sammler F aufgetragen, siehe Figur 1.

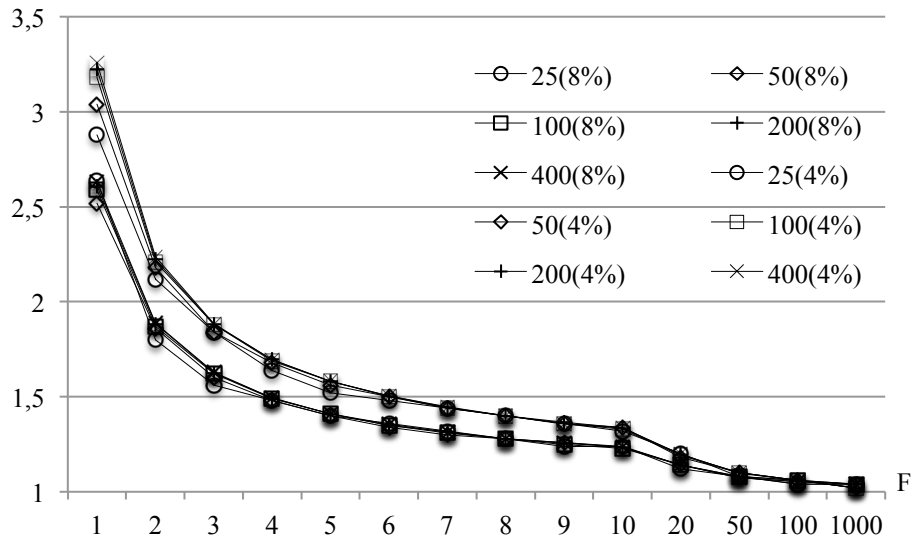

**Fig. 1.** Faktor f in Abhängigkeit der Anzahl Sammler F für verschiedene Albumgrößen B und Nachkaufprozentsätze K/B

Der Faktor f wird wie erwartet mit der Anzahl Sammler langsam kleiner. Allerdings strebt er für das WM-Album für großes F nicht wie nach (1) erwartet gegen ln(ln(640))≈1.87, sondern gegen Eins!

Bei der Simulation sind die Faktoren aller Alben bei festem Nachkaufprozentsatz fast identisch. Das ist überraschend. Zwar ist der Prozentsatz der Nachkaufkarten fast gleich, aber im größten Album befinden sich 16-mal so viele Bilder wie im kleinsten. Nur bei wenigen Tauschpartnern gibt es geringe Unterschiede. Je mehr Tauschpartner da sind, umso besser ist die Approximation. Damit kann der Sammler also die Fairness eines Albums alleine mit diesem Prozentsatz einschätzen, ohne komplizierte Berechnungen oder Simulationen machen zu müssen.

Daher formulieren wir als Hauptergebnis unsere

**Vermutung 1:** Das Verhältnis der mittleren Anzahl zu kaufender Karten pro Sammler zur Albumgröße hängt in sehr guter Näherung nur vom Anteil der nachkaufbaren Karten ab, aber nicht von der Albumgröße.

Wir haben mit einigen auf diesem Gebiet aktiven Wissenschaftlern korrespondiert, um herauszufinden, ob dieses Ergebnis in der Literatur bekannt und vielleicht schon bewiesen ist [9]. Bisher war das Ergebnis negativ.

## 3  Ansätze zum Beweis der Vermutung

Wir fangen mit einem Einzelsammler (F=1) an. Ein Sammler muß ohne Nachkaufen im Mittel $BH(B)$ kaufen. Wenn er die letzten K Karten nachkaufen kann, spart er also insgesamt $BH(K)$ Karten, wobei wir die Kosten des Nachkaufs erst mal vernachlässigen. Dies bedeutet also, dass der Faktor durch Nachkaufen für großes B näherungsweise auf

$$f = H(B) - H(K) + \frac{K}{B} \approx \ln B - \ln K + \frac{K}{B} = \ln \frac{B}{K} + \frac{K}{B} \quad (2)$$

reduziert wird, d. h. für einen Sammler hängt der Faktor näherungsweise nur vom Prozentsatz der Nachkaufkarten zur Albumgröße ab. Allerdings sieht man auch durch Vergleich mit den exakten Werten [3], dass das Ergebnis für einen Sammler nicht exakt stimmt, d. h. dass die Faktoren wirklich unterschiedlich sind, das heißt, die Unterschiede kommen nicht nur durch die Simulation zustande. Aber es ist interessant, dass die Unterschiede so gering sind, obwohl wir nur 1000 Simulationen durchgeführt haben.

Aber Formel (2) liefert noch eine weitere Erkenntnis, und zwar kann man umgekehrt näherungsweise bestimmen, wie viele Bilder K man im Mittel nachkaufen müsste, wenn man M Bilder gekauft hat. Dazu vernachlässigt man in (2) die Nachkaufbilder, multipliziert mit B und löst nach K auf

$$K \approx B \exp\left(-\frac{M}{B}\right) \quad (3)$$

und weiß damit auch, dass man im Mittel etwa B-K Bilder im Album hat und M-(B-K) Doppelte hat.

Wir wollen jetzt noch prüfen, ob es für die Standardabweichung einen ähnlichen Zusammenhang gibt. Durch einfaches Einsetzen erhalten wir für die Varianz für einen Sammler ohne Nachkaufen

$$V = \sum_{i=1}^{B} \frac{(B-i)/B}{i^2/B^2} = B^2 \sum_{i=1}^{B} \frac{1}{i^2} - B \sum_{i=1}^{B} \frac{1}{i} = B^2 H_2(B) - BH(B)$$

wobei $H_2$ die verallgemeinerte harmonische Zahl zweiter Ordnung ist. Wir versuchen daher, für die Varianz des Einzelsammlers mit Nachkaufen dasselbe Argument wie bei dem Erwartungswert anzuwenden, nämlich dass die letzten K Terme in der Summe der Varianzen der Wartezeiten wegfallen. Damit erhalten wir

$$V = B^2\left(H_2(B) - H_2(K)\right) - B\left(H(B) - H(K)\right)$$

und wir machen eine ähnliche Approximation wie in (2) und ersetzen die Summen durch endliche Integrale, insbesondere

$$H_2(n) = \sum_{i=1}^{n} \frac{1}{i^2} \approx \int_{1}^{n} \frac{dx}{x^2} = 1 - \frac{1}{n}$$

Für die Abschätzung der Standardabweichung der Anzahl der zu sammelnden Bilder kann man die Varianz durch B normieren und erhält schließlich

$$\frac{\sigma}{\sqrt{B}} \approx \sqrt{\left(\frac{B}{K} - 1\right) - \ln\left(\frac{B}{K}\right)}$$

Dadurch erkennt man, dass die normierte Standardabweichung wieder nur von dem Anteil K/B der Nachkaufbilder abhängig ist! Wir können also aufstellen:

**Vermutung 2:** Die normierte Standardabweichung hängt näherungsweise nur vom Prozentsatz der Nachkaufbilder K/B und der Anzahl der Tauschpartner F ab und nicht von der Größe B des Sammelalbums.

Wir haben uns nun gefragt, was passieren würde, wenn man ein Display mit D Bildern kaufen könnte, in dem wirklich keine Doppelten vorkommen, wie in [2] für die Schweizer WM-Bilder behauptet wurde. In diesem Fall spart man weniger als durch das Nachkaufen, aber man erhält die ersten D Bilder ohne Doppelte, sodass sich für den Faktor folgendes ergibt

$$f = H(B) - H(K) - (H(B) - H(B-D)) \approx \ln\left(\frac{B}{K}\right) - \ln\left(\frac{B}{B-D}\right) = \ln\left(\frac{B-D}{K}\right)$$

Dies ist kein völlig unerwartetes Ergebnis, denn es bedeutet, dass der Faktor nur vom Verhältnis der nach dem Display noch zu sammelnden Bilder B-D und den Nachkaufbildern K abhängt. Durch die D Bilder ohne Doppelte wird einfach die Albumgröße reduziert. Daher formulieren wir

**Vermutung 3:** Im Fall, dass man die ersten D Bilder ohne Doppelte erhalten kann, hängt der Faktor f näherungsweise nur vom Verhältnis der restlichen zu sammelnden Bilder B-D zu den Nachkaufbildern K und der Anzahl der Tauschpartner F ab und nicht von der Größe B des Sammelalbums.

Beide Vermutungen konnten wir mit Simulationen auch für mehrere Sammler bestätigen.

Für mehrere Sammler gestaltet sich der Nachweis wesentlich schwieriger, denn es ist im allgemeinen nicht bekannt, um wieviel sich die mittlere Anzahl der zu kaufenden Bilder durch das Nachkaufen reduziert. Wir können nur aus den Ergebnissen unserer

Simulationen und Formel (1) schließen, dass der Einfluß erheblich ist, da der Faktor f gegen Eins statt gegen ln(ln(B)) nach Formel (1) geht. Der Ansatz, mit dem Formel (1) bewiesen wird, ist auch wesentlich komplizierter [4] und wir haben bisher keinen Ansatz für einen Beweis unserer Vermutungen für mehrere Sammler gefunden.

## 4     Einfache Berechnung der mittleren Kosten und Streuung

Unsere Vernutung erlaubt eine einfache Abschätzung der mittleren Kosten eines Albums nur mit Hilfe der Grundrechenarten. Denn da der Faktor zumindest näherungsweise unabhängig von der Albumgröße ist, kann man den Zusammenhang zwischen der Anzahl der Sammler, dem Prozentsatz und dem Faktor einfach tabellieren bzw. in einem einzigen Diagramm darstellen und ablesen. In Tabelle 1 haben wir ein solches Diagramm abgebildet, das für die Albumgröße 100 mit jeweils 40.000 Simulationen erstellt wurde. Bei dieser Anzahl von Simulationen erhält man eine Genauigkeit von 0.01 auf einem Konfidenzniveau von 95% [1]. Eine höhere Genauigkeit ist nicht notwendig, da die Standardabweichung relativ groß ist.

| Sammler F | Prozentsatz der Nachkaufbilder K/B | | | | | |
|---|---|---|---|---|---|---|
| | 4% | 8% | 12% | 16% | 20% | 24% |
| 1 | 3,18 | 2,59 | 2,19 | 1,96 | 1,79 | 1,65 |
| 2 | 2,19 | 1,84 | 1,61 | 1,49 | 1,39 | 1,31 |
| 3 | 1,85 | 1,58 | 1,42 | 1,33 | 1,25 | 1,2 |
| 4 | 1,67 | 1,45 | 1,32 | 1,25 | 1,18 | 1,14 |
| 5 | 1,56 | 1,37 | 1,26 | 1,2 | 1,14 | 1,11 |
| 10 | 1,33 | 1,23 | 1,13 | 1,09 | 1,06 | 1,04 |
| 100 | 1,06 | 1,05 | 1,01 | 1 | 1 | 1 |

**Tab. 1.** Faktor f in Abhängigkeit der Anzahl Sammler F und verschiedenen Nachkaufprozentsätzen K/B

Eine analoge Tabelle können wir auf Grund von Vermutung 2 für die normierte Standardabweichung aufstellen, siehe Tabelle 2. Hier haben wir allerdings nicht mit einer festen Genauigkeit für die Schaätzung der Standardabweichung gearbeitet, sondern die Ergebnisse aus den Simulationen für Tabelle 1 abgeleitet. Dies bedeutet, dass die Genauigkeit geringer sein kann.

|  Sammler F | Prozentsatz der Nachkaufbilder K/B | | | | | |
|---|---|---|---|---|---|---|
|    | 4%   | 8%   | 12%  | 16%  | 20%  | 24%  |
| 1  | 4,56 | 2,99 | 2,28 | 1,79 | 1,53 | 1,3  |
| 2  | 4,2  | 2,75 | 2,1  | 1,7  | 1,43 | 1,21 |
| 3  | 4    | 2,73 | 2,06 | 1,67 | 1,38 | 1,19 |
| 4  | 4,12 | 2,73 | 2,06 | 1,66 | 1,37 | 1,17 |
| 5  | 4,2  | 2,73 | 2,08 | 1,67 | 1,37 | 1,16 |
| 10 | 4,34 | 2,87 | 2,17 | 1,69 | 1,36 | 1,08 |
| 100| 6,2  | 3,4  | 1,95 | 1,1  | 0,57 | 0,25 |

**Tab. 2.** Normierte Standardabweichung f in Abhängigkeit der Anzahl Sammler F und verschiedenen Nachkaufprozentsätzen K/B

Zusätzlich können wir einen weiteren Effekt berücksichtigen, nämlich dass die F Sammler eine Anzahl N von Nichtsammlern überreden könnten, für sie nachzukaufen. Dies kann praktisch für die Sammler sehr vorteilhaft sein. Da wir faires Tauschen vorausgesetzt haben, nehmen wir auch hier an, dass der Vorteil gerecht geteilt wird, d. h. der Prozentsatz $K/B$ wird noch einmal durch den Zusatzfaktor $(F+N)/F$ vergrößert.

Um die mittleren Kosten für ein Album zu bestimmen, muss man nun den Faktor f mittels Tabelle 1 bestimmen, bei abweichenden Prozentsätzen z. B. durch Interpolation zwischen den Kurven. Dabei muss man berücksichtigen, dass bei uns im Faktor f die Nachkaufkarten einfach berücksichtigt wurden. Dann kann man die Kosten für das Album (ohne Berücksichtigung des Nachkaufpreises) berechnen. Es fehlt jetzt noch der Aufpreis für die Nachkaufkarten, der hinzu addiert werden muss, insgesamt ergibt sich

$$Bf\frac{p}{P} + \frac{F+N}{F}K\left(b - \frac{p}{P}\right) \qquad (4)$$

Ein Beispiel soll diese demonstrieren: Für das Panini-WM-Album nehmen wir an, dass F=3 Freunde gemeinsam sammeln. Die Fairness beträgt etwa 8%, aus Tabelle 1 ermitteln wir f=1.58. Damit kostet das Album in guter Übereinstimmung mit [2] im Mittel etwa 123€. Aus Tabelle 2 können wir für die zugehörige normierte Standardabweichung 2,73 ermitteln, durch Multiplikation mit √640 erhalten wir eine Standardabweichung von etwa 69 Bildern, dies entspricht etwa 8,29€. D. h. mit etwa 95% Wahrscheinlichkeit muss jeder der 3 Sammler zwischen 107€ und 139€ ausgeben, um sein Album zu komplettieren. Überreden die Freunde jetzt einen Nichtsammler (N=1), für sie 50 Extrakarten nachzukaufen, so erhöht sich der Prozentsatz auf 10.6%, durch Interpolation ermitteln wir f=1.47. Wir müssen jetzt aber berücksichtigen, dass jeder Sammler etwa K=67 Karten nachkaufen (und bezahlen) muss, und erhalten etwa 117€ als mittlere Kosten für das Album. Hätte der Anbieter wie einige andere einen Nachkaufprozentsatz von 16% angeboten, d. h. 100 Bilder pro Sammler, so wäre f=1.32, so hätte das Album nur 107€ gekostet (N=0).

## 5    Zusammenfassung

Wir haben in dieser Arbeit mehrere neue Vermutungen aufgestellt und begründet, die es Sammlern erlauben, eine einfache Abschätzung der mittleren Kosten eines Sammelalbums sowie deren Standardabweichung nur mit Grundrechenarten (vgl. Formel 4) und zwei Tabellen durchzuführen. Durch Beispiele wurde die Nützlichkeit nachgewiesen. Nach unserer Recherche ist dies der erste veröffentlichte Ansatz für die Berechnung der Kosten, bisher wurden nur Simulationen [2] verwendet.

Wir konnten die Korrektheit unserer Vermutungen für den Spezialfall eines Sammlers nachweisen, aber nicht für den allgemeinen Fall. Im allgemeinen Fall haben wir anhand von umfangreichen Simulation die Plausibilität unserer Vermutungen nachgewiesen.

Damit ergibt sich auch eine interessante Betrachtung der Asymptotik des Faktors f im Sammelbilderproblem, siehe Tabelle 2. Ohne Nachkaufen wächst der Faktor mit der Albumgröße B gegen Unendlich. Mit Nachkaufen bleibt der Faktor für einen Einzelsammler bei festem Nachkaufprozentsatz dagegen beschränkt, selbst wenn B wächst. Und bei fairem Tauschen mit Nachkaufen strebt der Faktor gegen Eins, wenn F gegen Unendlich strebt, d. h. aus dem fairen Tauschen wird das ideale Tauschen, bei dem man sofort jedes Doppelte tauschen kann.

| Sammelbilderproblem | Ohne Nachkaufen | Mit Nachkaufen |
|---|---|---|
| **Einzelsammler (F=1)** | $\ln B$<br>Klassische Ergebnisse [3] | $\ln \frac{B}{K} + \frac{K}{B}$<br>siehe (2) |
| **Faires Tauschen (F>1)** | $\ln \ln B$<br>für großes F [4] | $\approx 1$<br>für großes F, siehe Figur 1 und Tabelle 2 |

**Tab. 3.** Asymptotik des Faktors bei verschiedenen Varianten des Sammelbilderproblems

Dies zeigt, dass im Vergleich zwischen Nachkaufen und Tauschen zumindest asymptotisch das Nachkaufen einen stärkeren Effekt hat als das Tauschen.

cherche, insbesondere der englischsprachigen Literatur, der Kommunikation mit Wissenschaftlern sowie der Formatierung dieses Texts unterstützt.

## Referenzen